\title{Limits of randomly grown graph sequences}
\author{Christian Borgs\footnote{Microsoft Research, Cambridge, MA}, Jennifer Chayes${}^*$, L\'aszl\'o Lov\'asz\footnote{Institute of Mathematics,
E\"otv\"os Lor\'and University, Budapest, Hungary. Research sponsored by OTKA Grant No.~67867.}~,\\
Vera S\'os\footnote{A.~R\'enyi Institute of Mathematics, Budapest,
Hungary}, Katalin Vesztergombi${}^\dag$}
\date{February 2009}

\documentclass{article}
\usepackage{amssymb,amsmath,hyperref,graphicx,theorem,delarray}

\nonstopmode

\long\def\killtext#1{}
\newtheorem{theorem}{Theorem}[section]
\newtheorem{prop}[theorem]{Proposition}
\newtheorem{lemma}[theorem]{Lemma}
\newtheorem{claim}[theorem]{Claim}

\theorembodyfont{\rmfamily}
\newtheorem{example}{Example}

\newtheorem{remark}[theorem]{Remark}

\newenvironment{proof}{\medskip\noindent{\bf Proof. }}{\hfill$\square$\medskip}

\begin{document}

\addtolength{\baselineskip}{3pt} \setlength{\oddsidemargin}{0.2in}

\def\R{{\mathbb R}}
\def\one{{\mathbf 1}}
\def\Q{{\mathbb Q}}
\def\Z{{\mathbb Z}}
\def\C{{\mathbb C}}
\def\hom{{\rm hom}}
\def\inj{{\rm inj}}
\def\surj{{\rm sur}}
\def\Inj{{\rm Inj}}
\def\neighb{{\rm neighb}}
\def\inj{{\rm inj}}
\def\ind{{\rm ind}}
\def\sur{{\rm sur}}
\def\PAG{{\rm PAG}}
\def\SPAG{{\rm SPAG}}
\def\eps{\varepsilon}

\def\Ge{{\mathbb G}}
\def\Ha{{\mathbb H}}

\def\Ab{{\mathbf A}}
\def\Sb{{\mathbf S}}

\def\tv{{\rm tv}}
\def\tr{{\rm tr}}
\def\cost{\hbox{\rm cost}}
\def\val{\hbox{\rm val}}
\def\rk{{\rm rk}}
\def\diam{{\rm diam}}
\def\Ker{{\rm Ker}}
\def\QG{{\cal QG}}
\def\QGM{{\cal QGM}}
\def\CD{{\cal CD}}
\def\Pr{{\sf P}}
\def\E{{\sf E}}
\def\Var{{\sf Var}}
\def\Ent{{\sf Ent}}
\def\T{{^\top}}

\def\AA{{\cal A}}\def\BB{{\cal B}}\def\CC{{\cal C}}
\def\DD{{\cal D}}\def\EE{{\cal E}}\def\FF{{\cal F}}
\def\GG{{\cal G}}\def\HH{{\cal H}}\def\II{{\cal I}}
\def\JJ{{\cal J}}\def\KK{{\cal K}}\def\LL{{\cal L}}
\def\MM{{\cal M}}\def\NN{{\cal N}}\def\OO{{\cal O}}
\def\PP{{\cal P}}\def\QQ{{\cal Q}}\def\RR{{\cal R}}
\def\SS{{\cal S}}\def\TT{{\cal T}}\def\UU{{\cal U}}
\def\VV{{\cal V}}\def\WW{{\cal W}}\def\XX{{\cal X}}
\def\YY{{\cal Y}}\def\ZZ{{\cal Z}}

\maketitle

\tableofcontents

\begin{abstract}
Motivated in part by various sequences of graphs growing under random
rules (like internet models), convergent sequences of dense graphs
and their limits were introduced by Borgs, Chayes, Lov\'asz, S\'os
and Vesztergombi and by Lov\'asz and Szegedy. In this paper we use
this framework to study one of the motivating class of examples,
namely randomly growing graphs. We prove the (almost sure)
convergence of several such randomly growing graph sequences, and
determine their limit. The analysis is not always straightforward: in
some cases the cut distance from a limit object can be directly
estimated, in other case densities of subgraphs can be shown to
converge.
\end{abstract}

\section{Introduction}

Convergent graph sequences and their limits have been studied in
connection with internet models, statistical physics, extremal graph
theory, and more. In the context of dense graphs, a rather complete
theory has emerged. One can define a notion of convergence based on
the convergence of densities of subgraphs. An appropriate notion of
distance between two graphs, called their {\it cut distance}, can be
defined, so that convergent sequences are Cauchy in this metric and
vice versa. The completion of the metric space of graphs relative to
this metric can be described, and its elements, i.e., limit objects
for convergent graph sequences, can be characterized in various ways.
To mention one of these, limit objects can be described by 2-variable
symmetric measurable functions $[0,1]^2\to[0,1]$.

The goal of this paper is study in this framework one of the
motivating class of examples, namely randomly growing graphs.
Typically, such a sequence of graphs grows by every now and then
adding a new node, and then creating new edges (between the new node
and the old ones, or between two old nodes) randomly, from some
simple distribution determined by local conditions.

We will prove the (almost sure) convergence of several such randomly
growing graph sequences, and determine their limit. This analysis is
not always straightforward: in some cases the cut distance from a
limit object can be directly estimated, in other case densities of
subgraphs can be shown to converge.

\section{Preliminaries}

In this section we summarize those notions and results concerning
convergent graph sequences and their limits which are relevant for
the rest of the paper.

\subsection{Convergent graph sequences}

For two simple graphs $F$ and $G$, $\hom(F,G)$ denotes the number of
homomorphisms (adjacency-preserving maps) from $V(F)$ to $V(G)$. We
also consider the {\it homomorphism densities}
\begin{equation}\label{HOMDENSE}
t(F,G)=\frac{\hom(F,G)}{|V(G)|^{|V(F)|}}.
\end{equation}
(Thus $t(F,G)$ is the probability that a random map of $V(F)\to V(G)$
is a homomorphism.)

A sequence $(G_n)$ of graphs is {\it convergent}, if the sequence
$t(F,G_n)$ has a limit for every simple graph $F$.

Convergent graph sequences have a limit object, which can be
represented as measurable functions \cite{LSz1}. Let $\WW$ denote the
space of all bounded measurable functions $W:~[0,1]^2\to \R$ such
that $W(x,y)=W(y,x)$ for all $x,y\in[0,1]$. We also define
$\WW_0=\{W\in \WW:~0\le W\le 1\}$. For every simple graph $F$ and
$W\in\WW$, we define
\[
t(F,W)=\int_{[0,1]^{V(F)}} \prod_{ij\in E(F)} W(x_i,x_j)\,dx.
\]

Every finite simple graph $G$ can be represented by a function
$W_G\in\WW_0$: Let $V(G)=\{1,\dots,n\}$. Split the interval $[0,1]$
into $n$ equal intervals $J_1,\dots,J_n$, and for $x\in J_i, y\in
J_j$ define
\[
W_G(x,y)=
  \begin{cases}
    1, & \text{if $ij\in E(G)$}, \\
    0, & \text{otherwise}.
  \end{cases}
\]
Informally, we replace the $(i,j)$ entry in the adjacency matrix of
$G$ by a square of size $(1/n)\times(1/n)$, and define the value of
the function $W_G$ on this square as the corresponding entry of the
adjacency matrix.

Graphons represent limits of convergent graph sequences in the
following sense.

\begin{theorem}\label{THM:LIMIT}
{\rm(a)} For every convergent graph sequence $(G_n)$ there is a
$W\in\WW$ such that $t(F,G_n)\to t(F,W)$ for every simple graph $F$.

\smallskip

{\rm(b)} This function $W$ is uniquely determined up to measure
preserving transformations in the following sense: For every other
limit function $W'$ there are measure preserving maps
$\phi,\psi:~[0,1]\to[0,1]$ such that
$W(\phi(x),\phi(y))=W'(\psi(x),\psi(y))$.

\smallskip

{\rm(c)} Every function $W\in\WW_0$ arises as the limit of a
convergent graph sequence.
\end{theorem}

Parts (a) and (c) of the theorem were proved in \cite{LSz1}, and part
(b), in \cite{BCL}. The proof of (c) in \cite{LSz1} depends on
$W$-random graphs, to be discussed in the next section.

We could consider any probability space $(\Omega,\AA,\pi)$ instead of
$[0,1]$, with a symmetric measurable function
$W:~\Omega\times\Omega\to [0,1]$. These structures are called {\it
graphons}. The densities $t(F,W)$ in a graphon could be defined by a
similar integral. Considering graphons would not give greater
generality, since we could always replace $(\Omega,\AA,\pi)$ by the
uniform measure on $[0,1]$. Still, it is sometimes useful to
represent the limit object by other probability spaces, as we shall
see.

\subsection{Distance of graphs}

The {\it cut-norm} introduced in \cite{FK} is defined for $W\in\WW$
by
\[
\|W\|_\square=\sup_{S,T\subset [0,1]} \Bigl|\int_{S\times T}
W(x,y)dxdy\Bigr|,
\]
where the supremum goes over measurable subsets of $[0,1]$. We define
the {\it cut-distance} of two functions in $\WW$ by
\begin{equation}
\label{d-cut-graphon}
\delta_\square(U,W)
=\inf_{\phi:\,[0,1]\to[0,1]} \|U-W^\phi\|_\square
\end{equation}
where the infimum goes over all invertible maps $\phi:~[0,1]\to[0,1]$
such that both $\phi$ and its inverse are measure preserving, and
$W^\phi$ is defined by $W^\phi(x,y)=W(\phi(x),\phi(y))$. For two
graphs $G$ and $G'$, this yields a distance
\[
\delta_\square(G,G')=\delta_\square(W_G,W_{G'}).
\]

\begin{remark}\label{RM:DIST}
(a) We call this a ``distance'' rather than a ``metric'' since two
different graphs can have distance $0$. This is the case when one
graph can be obtained from the other by replacing each node by the
same number of twins, or more generally, when both can be obtained
from a third graph this way. To get a metric, we should identify such
pairs of graphs. Similarly, to get a metric on $\WW_0$, we have to
identify functions $U,W$ for which $\delta_\square(U,W)=0$. Several
characterizations of such pairs are given in \cite{BCL}.

\smallskip

(b) There are combinatorial, but somewhat lengthy ways to define this
distance between graphs; see \cite{BCLSV1}.
\end{remark}

We can define a similar distance function based on other norms. We
shall use the $L_1$-norm
\[
\|W\|_1=\int_{[0,1]^2} |W(x,y)|\,dx\,dy,
\]
from which we can define the {\it edit distance} of two functions in
$\WW$ by
\begin{equation}
\label{d-edit-graphon}
\delta_1(U,W)
=\inf_{\phi:\,[0,1]\to[0,1]} \|U-W^\phi\|_1
\end{equation}

The following characterization of convergent graph sequences was
proved in \cite{BCLSV1} (see \cite{BCLSV2} for other
characterizations not used in this paper).

\begin{theorem}\label{THM:CONV-DIST}
A sequence of graphs $(G_n)$ is convergent if and only if it is
Cauchy in the $\delta_\square$ distance. The sequence $(G_n)$
converges to $W$ if and only if $\delta_\square(W_{G_n},W)\to 0$.
Furthermore, there is a way to label the nodes of the graphs in the
sequence so that $\|W_{G_n}-W\|_\square\to 0$.
\end{theorem}

If the graphs $G_n$ are labeled so that $\|W_{G_n}-W\|_\square\to 0$,
then
\[
\sup_{S,T}\left|\int_{S\times T} (W_{G_n}-W)\right| \to 0\qquad
(n\to\infty)
\]
In particular, it follows that
\begin{equation}\label{EQ:STW}
\int_{S\times T} (W_{G_n}-W) \to 0
\end{equation}
for every product set $S\times T$, which implies that $W_{G_n}\to W$
in the weak* topology of $L_\infty([0,1]^2)$. Convergence in the norm
$\|.\|_\square$ is, however, not equivalent to convergence in this
weak* topology, as the sequence prefix attachment graphs shows
(Section \ref{SEC:PFX-ATT}).

\subsection{$W$-random graphs and extensions}\label{SEC:WRAND}

Let $(\Omega,\AA,\pi,W)$ be a graphon. For every finite subset
$S\subseteq \Omega$ we define two graphs $\Ge(S,W)$ and $\Ha(S,W)$ on
$V(\Ge(S,W))=V(\Ha(S,W))=S$. In $\Ge(S,W)$, we connect $i,j\in S$,
$i\not= j$ with probability $W(i,j)$. In $\Ha(S,W)$, we connect
$i,j\in S$, $i\not= j$ by an edge with weight $W(i,j)$. If $W$ is
$\{0,1\}$-valued, then $\Ge(S,W)=\Ha(S,W)$ is deterministic, and can
be considered as an ``induced subgraph''.

Let $\Sb_n$ be a random $n$-element subset of $\Omega$ (each element
of $\Sb_n$ chosen independently from the distribution $\pi$). The
graph $\Ge(n,W)=\Ge(\Sb_n,W)$ is called a {\it $W$-random} graph. The
following fact was shown in \cite{LSz1} (for the case when the
underlying probability space is the uniform distribution on $[0,1]$,
but this is no essential restriction of generality).

\begin{lemma}\label{LEM:WRAND}
With probability $1$, the sequence $\Ge(n,W)$ is convergent and its
limit is represented by the function $W$.
\end{lemma}

\killtext{
It will be convenient to use the following notation.
Consider a graph $F$ whose edges are partitioned into two sets $E'$
and $E''$, called ``blue'' and ``red'', to get a two-edge-colored
graph $F'=(V,E',E'')$. We define
\begin{equation} \label{tindFW-def}
t(F',W)=\int\limits_{\Omega^V} \prod_{ij\in E'}
W(x_i,x_j)\prod_{ij\in E''} (1-W(x_i,x_j))\,\prod_{i\in
V}d\pi(x_i)\,.
\end{equation}
If all edges are blue, then $t(F',W)=t(F,W)$. If all edges are red,
then $t(F',W)=t(F,1-W)$. In general, $t(F',W)$ can be expressed as
\begin{equation}\label{SIEVE}
t(F',W) = \sum_{Y\subseteq E''} (-1)^{|E''\setminus Y|} t(F\setminus
Y, W).
\end{equation}
For a simple graph $F$, let $K_F$ denote the complete graph on
$V(F)$, where the edges in $E(F)$ are colored red, the other edges
are colored blue. Then for any two simple graphs $F$ and $G$,
$t(K_F,W_G)$ is the probability that a random map $V(F)\to V(G)$
preserves both adjacency and non-adjacency. With this notation, we
can express the probability that $\Ge(n,W)$ is a given graph $F$
(with a fixed labeling of the nodes):
\[
\Pr(\Ge(n,W)=F)=t(K_F,W).
\]
}

In this paper, we will also need sequences $S_n$ of subsets of
$\Omega$ that are not random, but still $\Ge(S_n,W)$ converges to
$W$. We prove and use the following sufficient condition for a
deterministic sequence $S_n$. Let $(\Omega,d)$ be a metric space, and
$\pi$, a probability measure on the Borel subsets of $(\Omega,d)$.
For every $n\ge 1$, let $S_n$ be a finite subset of $\Omega$ such
that $|S_n|\to \infty$. We say that the sequence $(S_n)$ is {\it well
distributed in a set} $X\subseteq\Omega$, if $|S_n\cap X|/|S_n|\to
\pi(X)$ as $n\to\infty$. We say that $(S_n)$ is {\it well distributed
in $(\Omega,d,\pi)$}, if for every $\eps>0$ there exists a partition
$\{P_1,\dots,P_m\}$ of $\Omega$ into sets with diameter at most
$\eps$ such that $S_n$ is well distributed in each $P_j$.

\begin{lemma}\label{LEM:CONV}
Let $(\Omega,d,\pi)$ be a metric space with an atom-free probability
measure. Let $W:~\Omega\times\Omega \to[0,1]$ be a symmetric
measurable function that is almost everywhere continuous. Let $S_n$
be a sequence of sets that is well distributed in $(\Omega,d,\pi)$.

\smallskip

{\rm (a)} Then $\delta_1(W_{\Ha(S_n,W)},W)\to 0$ and with probability
$1$, $\delta_\square(W_{\Ge(S_n,W)},W)\to 0$.

\smallskip

{\rm (b)} If $W$ is $0$-$1$ valued, then
$\delta_1(W_{\Ge(S_n,W)},W)\to 0$.
\end{lemma}

\noindent It is clear that such a conclusion cannot hold without some
assumption on $W$, since a general measurable function could be
changed on the sets $S_n\times S_n$ arbitrarily without changing its
subgraph densities.

\begin{proof}
(a) First we construct a special partition of $\Omega$.

\begin{claim}\label{CLAIM:EQUI}
There exists a sequence of partitions $\QQ_n$ of $\Omega$ into
$|S_n|$ sets such that every partition class contains exactly one
point of $S_n$, the maximum diameter of partition classes tends to
$0$, and the maximum of $\bigl|\pi(Q)|S_n|-1\bigr|$ $(Q\in\QQ_n)$,
tends to $0$.
\end{claim}

Let $\eps>0$. Consider a partition $\{P_1,\dots,P_m\}$ into sets with
diameter at most $\eps$ such that $S_n$ is well distributed in every
$P_j$. For $n$ large enough, we have $(1-\eps)\pi(P_j) \le |S_n\cap
P_j|/|S_n|\le (1+\eps)\pi(P_j)$ for every $j$. Let us partition each
set $P_j$ into $|S_n\cap P_j|$ sets of equal measure, each containing
exactly one point of $S_n\cap P_j$ to get the partition $\QQ_n$. It
is clear that this sequence of partitions has the properties as
required in the Claim.

For each $n$ and $s\in S_n$, let $Q_s$ be the partition class of
$\QQ$ containing $s$. Define the function $W_n$ as follows: for
$s,s'\in S_n$ and $(x,y)\in Q_s\times Q_{s'}$, let
$W_n(x,y)=W(s,s')$. Then $W_n(x,y)\to W(x,y)$ in every point $(x,y)$
where $W$ is continuous, in particular $W_n\to W$ almost everywhere.
This implies that
\begin{equation}\label{WWN}
\|W_n-W\|_1\to 0\qquad(n\to\infty).
\end{equation}

We can view $W_n$ as $W_{H_n}$, where $H_n$ is a weighted graph with
$V(H_n)=S_n$, the weight of node $s\in S_n$ is $\pi(Q_s)$, and the
weight of  $ss'$ ($s,s'\in S$) is $W(s,s')$. Note that $H_n$ is
almost the same weighted graph as $\Ha_n=\Ha(S_n,W)$: they are
defined on the same set of nodes, the edges have the same weights,
and the nodeweight $\pi(Q_s)$ is asymptotically $1/|S_n|$ by the
Claim. Given $\eps>0$, we have $|\pi(Q_s)-1/|S_n||<\eps/|S_n|$ if $n$
is large enough. Hence there is a measure preserving bijection
$\phi:~[0,1]\to[0,1]$ and a set $R\subseteq[0,1]$ of measure $\eps$
such that
\[
W_{H_n}(x,y)=W^\phi_{\Ha_n}(x,y) \qquad (x,y\notin R).
\]
This implies that
\begin{equation}\label{HNH}
\delta_1(H_n,\Ha_n)\to 0\qquad (n\to\infty).
\end{equation}
By Lemma 4.3 from \cite{BCLSV1} it follows that with probability $1$,
\begin{equation}\label{HNG}
\delta_\square(\Ha(S_n,W),\Ge(S_n,W))\to 0\qquad (n\to\infty).
\end{equation}
Equations \eqref{WWN}, \eqref{HNH} and \eqref{HNG} imply that
$\Ge(S_n,W)\to W$ with probability $1$.

(b) follows trivially, since in this case $\Ha(S_n,W)=\Ge(S_n,W)$.
\end{proof}

We note that (b) would also follow from the result of Pikhurko
\cite{Pik} that if a graph sequence tends to a 0-1 valued function
$W$ in the $\delta_\square$ distance, then it also tends to $W$ in
the $\delta_1$ distance.

\subsection{Pixel picture}

We have seen that every finite simple graph $G$ can be represented by
a function $W_G\in\WW_0$. In fact, this representation is very useful
for creating figures representing graphs.

\begin{figure}[htb]
  \centering
  \includegraphics*[width=.25\textwidth,angle=180,origin=c]{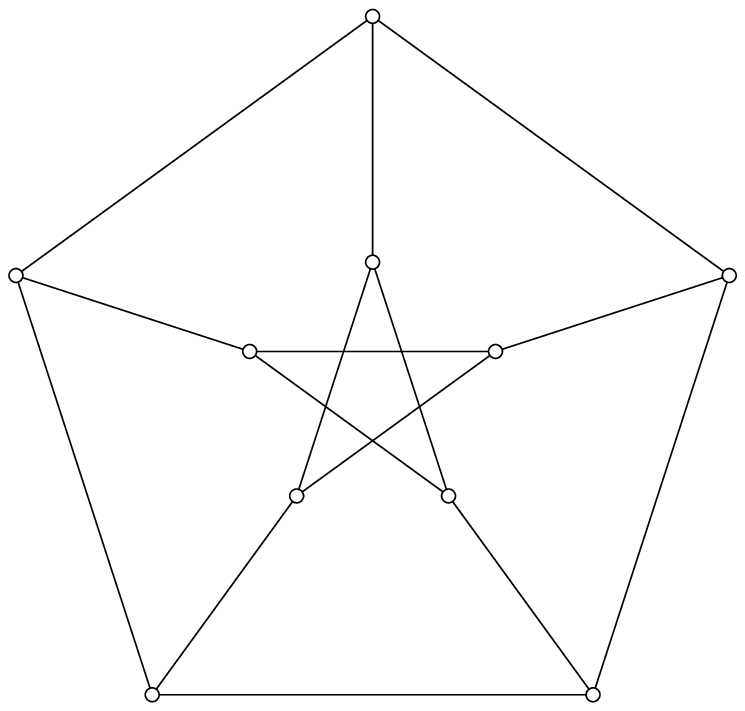}
  $\begin{array}[b]({cccccccccc})0&1&0&0&1&1&0&0&0&0\\
1&0&1&0&0&0&1&0&0&0\\
0&1&0&1&0&0&0&1&0&0\\
0&0&1&0&1&0&0&0&1&0\\
1&0&0&1&0&0&0&0&0&1\\
1&0&0&0&0&0&0&1&1&0\\
0&1&0&0&0&0&0&0&1&1\\
0&0&1&0&0&1&0&0&0&1\\
0&0&0&1&0&1&1&0&0&0\\
0&0&0&0&1&0&1&1&0&0
\end{array}$
  \includegraphics*[width=.25\textwidth]{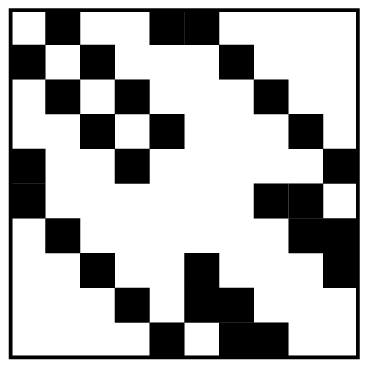}
  \caption{The Petersen graph, its adjacency matrix, and its pixel picture}
  \label{FIG:PIXEL}
\end{figure}

Every function $W\in\WW_0$ can be represented by a grayscale picture
on the unit square: the point $(x,y)$ is black if $W(x,y)=1$, it is
white if $W(x,y)=0$, and it is appropriately dark grey if
$0<W(x,y)<1$. For a graph, this picture gives a black-and-white
picture consisting of a finite number of ``pixels''. The origin is in
the upper left corner (as for a matrix). Figure \ref{FIG:PIXEL}
illustrates this construction. Note that the function associated with
a graph depends on the ordering of the nodes.

\begin{figure}[htb]
  \centering
  \includegraphics*[width=.3\textwidth]{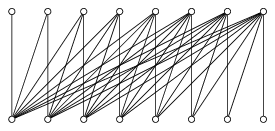}\qquad
  \includegraphics*[width=.25\textwidth]{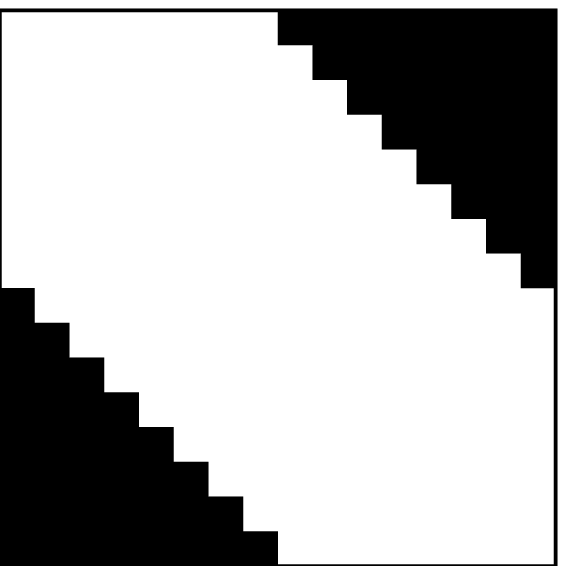}\qquad
  \includegraphics*[width=.25\textwidth]{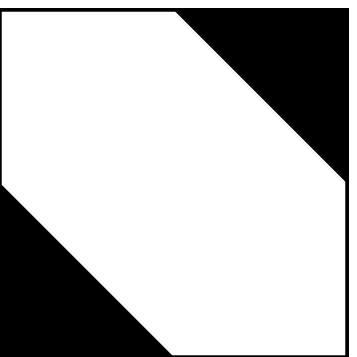}
  \caption{A half-graph, its pixel picture, and the limit function}
  \label{FIG:HALF}
\end{figure}

\begin{example}[Half graphs]\label{EXA:HALF}
Consider the {\it half-graphs} $H_{n,n}$: they are bipartite graphs
on $2n$ nodes $\{1,\dots,n,1',\dots,n'\}$, where $i$ is connected to
$j'$ if and only if $i\le j'$. It is easy to see that this sequence
is convergent, and to guess the limit function (Figure
\ref{FIG:HALF}).
\end{example}

\begin{example}[Erd\H{o}s-R\'enyi random graphs]\label{EXA:RAND}
The pixel picture of a random graph is essentially grey.
\end{example}

\begin{figure}[htb]
  \centering
  \includegraphics*[height=120pt,width=120pt]{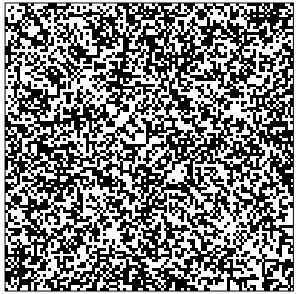}
  \caption{A random graph with 100 nodes and with edge density 1/2}
  \label{FIG:RANDOM}
\end{figure}

\begin{figure}[htb]
  \centering
  \includegraphics*[height=120pt,width=120pt]{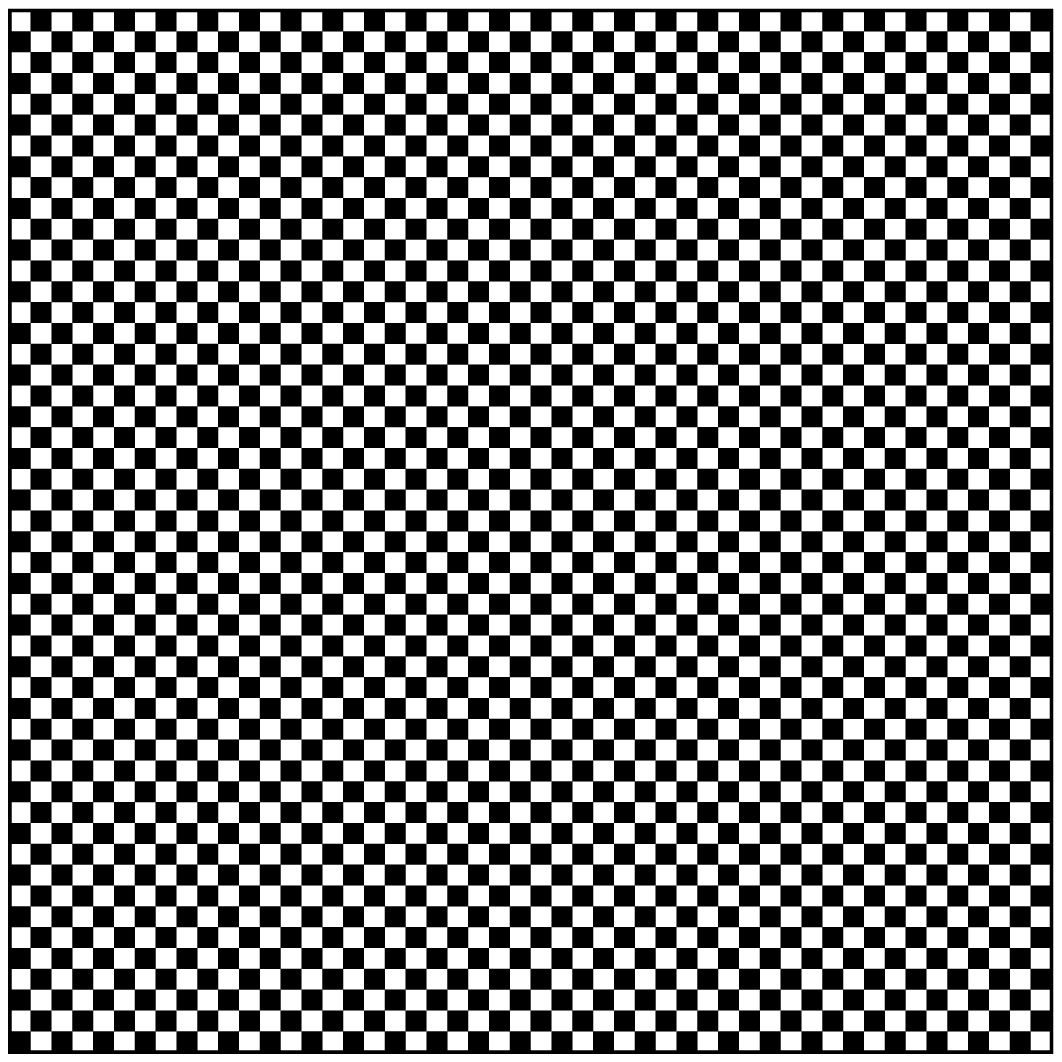}
  \includegraphics*[height=115pt,width=115pt]{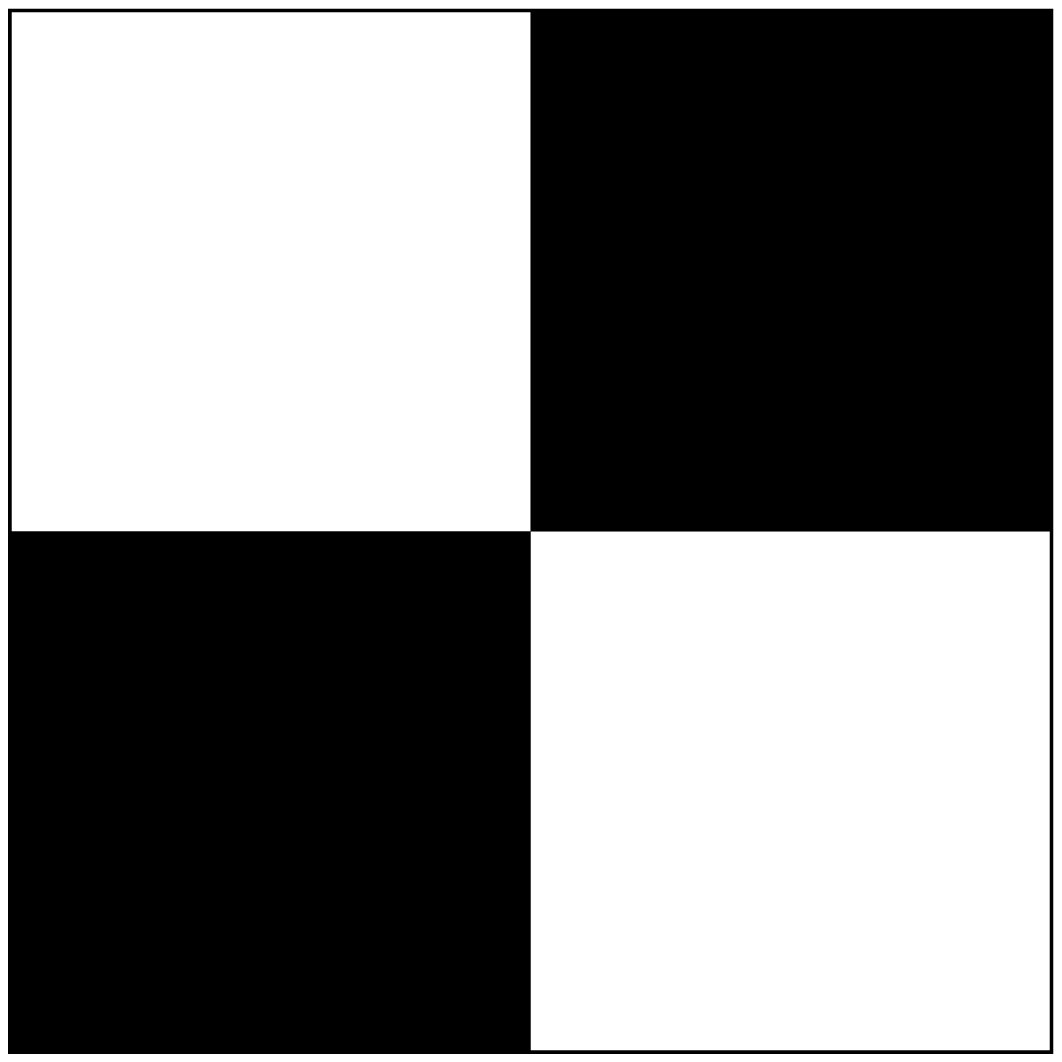}
  \caption{A chessboard and the pixel picture obtained by rearranging the
  rows and columns}\label{FIG:CHESS}
\end{figure}

The following simple example illustrates the importance of the
ordering of the nodes:

\begin{example}[Chessboard]\label{EXA:CHESS}
The $100\times 100$ chessboard in Figure \ref{FIG:CHESS} is the pixel
picture of a complete bipartite graph. It is also uniformly grey, so
one might assume that it represents a graph that is close to random.
But rearranging the rows and columns so that odd indexed columns come
first, we see that it is isomorphic to the graph represented by the
$2\times 2$ chessboard.

This example also shows that different graphs may be represented by
the same pixel picture: all complete bipartite graphs with equal
color classes have the same pixel picture. If we restrict our
attention to graphs with no twin nodes, the pixel picture will
determine the graph.
\end{example}

The pixel picture of a random graph remains uniformly grey, no matter
how you reorder the nodes.

It is easy to verify that
\[
t(F,G)=t(F,W_G)
\]
for every finite simple graph $G$.

\section{Convergent graph sequences and their limits}

\subsection{Growing uniform attachment graphs}\label{SEC:UNIF-ATT}

We generate a randomly growing graph sequence $G^{\rm ua}_n$ as
follows. We start with a single node. At the $n$-th iteration, a new
node is born, and then every pair of nonadjacent nodes is connected
with probability $1/n$. We call this graph sequence a {\it randomly
grown uniform attachment graph sequence}.

\begin{figure}[htb]
  \centering
  \includegraphics*[height=120pt,width=120pt]{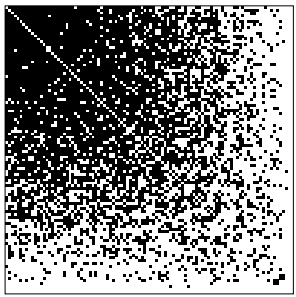}
  \caption{A randomly grown uniform attachment graph with 100 nodes}
  \label{FIG:UNIF}
\end{figure}

Let us do some simple calculations. After $n$ steps, let
$\{0,1,\dots,n-1\}$ be the nodes (born in this order). The
probability that nodes $i < j$ are not connected is
$\frac{j}{j+1}\cdot \frac{j+1}{j+2}\cdots\frac{n-1}{n}= \frac{j}{n}$.
These events are independent for all pairs $(i,j)$. The expected
degree of $j$ is
\[
\sum_{i=0}^{j-1}\frac{n-j}{n} + \sum_{i=j+1}^{n-1}\frac{n-i}{n} =
\frac{n-1}{2}-\frac{j(j-1)}{2n}.
\]
The expected number of edges is
\[
\frac{1}{2}\sum_{j=0}^{n-1}
\left(\frac{n-1}{2}-\frac{j(j-1)}{2n}\right) = \frac{n^2-1}{6}.
\]

To figure out the limit function, note that the probability that
nodes $i$ and $j$ are connected is $1-\max(i,j)/n$. If $i=xn$ and
$j=yn$, then this is $1-\max(x,y)$. This motivates the following:

\begin{theorem}\label{THM:MAX}
The sequence $G^{\rm ua}_n$ tends to the limit function $1-\max(x,y)$
with probability $1$.
\end{theorem}

\begin{proof}
For a fixed $n$, the events that nodes $i$ and $j$ are connected are
independent for different $i,j$, and so by the computation above,
$G^{\rm ua}_n$ has the same distribution as $\Ge(S_n,1-\max(x,y))$,
where $S_n=\{0, 1/n,\dots,(n-1)/n\}$. It is easy to see that this
sequence is well distributed in the metric space $[0,1]$ with the
uniform measure, and so the Theorem follows by Lemma \ref{LEM:CONV}.

One can get a good explicit bound on the convergence rate by
estimating the cut-distance of $W_{G^{\rm ua}_n}$ and $1-\max(x,y)$,
using the Chernoff-Hoeffding bound.
\end{proof}

\subsection{Growing ranked attachment graphs}\label{SEC:R-ATT}

This randomly growing graph sequence $G^{\rm ra}_n$ is generated
somewhat similarly. We start with a single node. At the $n$-th
iteration, a new node is born, and it is connected to node $i$ with
probability $1-i/n$. Then every pair of nonadjacent nodes is
connected with probability $2/n$. We call this graph sequence a {\it
randomly grown ranked attachment graph sequence}.

\begin{theorem}\label{THM:1-XY}
The sequence $G^{\rm ra}_n$ tends to the limit function $1-xy$ with
probability $1$.
\end{theorem}

\begin{proof}
The probability that nodes $i$ and $j$ are {\it not} connected after
the $n$-th step is
\begin{align*}
p_{ij}&=\frac{i}{j}\cdot \left(1-\frac2j\right)\cdot
\left(1-\frac2{j+1}\right)\cdots
\left(1-\frac2n\right)=\frac{i(j-2)(j-1)}{j(n-1)n}\\
&= \frac{ij}{n^2} - \frac{(3n-j)ij-2ni}{jn(n-1)} = \frac{ij}{n^2} -
q_{ij},
\end{align*}
where $0<q_{ij}<\min\{\frac3n,ij/n^2\}$. Furthermore, these events
are independent for different pairs $i,j$. Therefore, we can generate
the graph $G^{\rm ra}_n$ as follows: We generate $\Ge(S_n,1-xy)$,
where $S_n=\{0, 1/n,\dots,(n-1)/n\}$, and then connect each
nonadjacent $i$ and $j$ with probability $1-p_{ij}$. Since
$\Ge(S_n,1-xy)$ tends to the function $1-xy$ by Lemma \ref{LEM:CONV}
and the added edges change $\Ge(S_n,1-xy)$ negligibly in
$\delta_\square$ distance, the Theorem follows.
\end{proof}

\subsection{Growing prefix attachment graph}\label{SEC:PFX-ATT}

\begin{figure}[htb]
  \centering
  \includegraphics*[height=120pt]{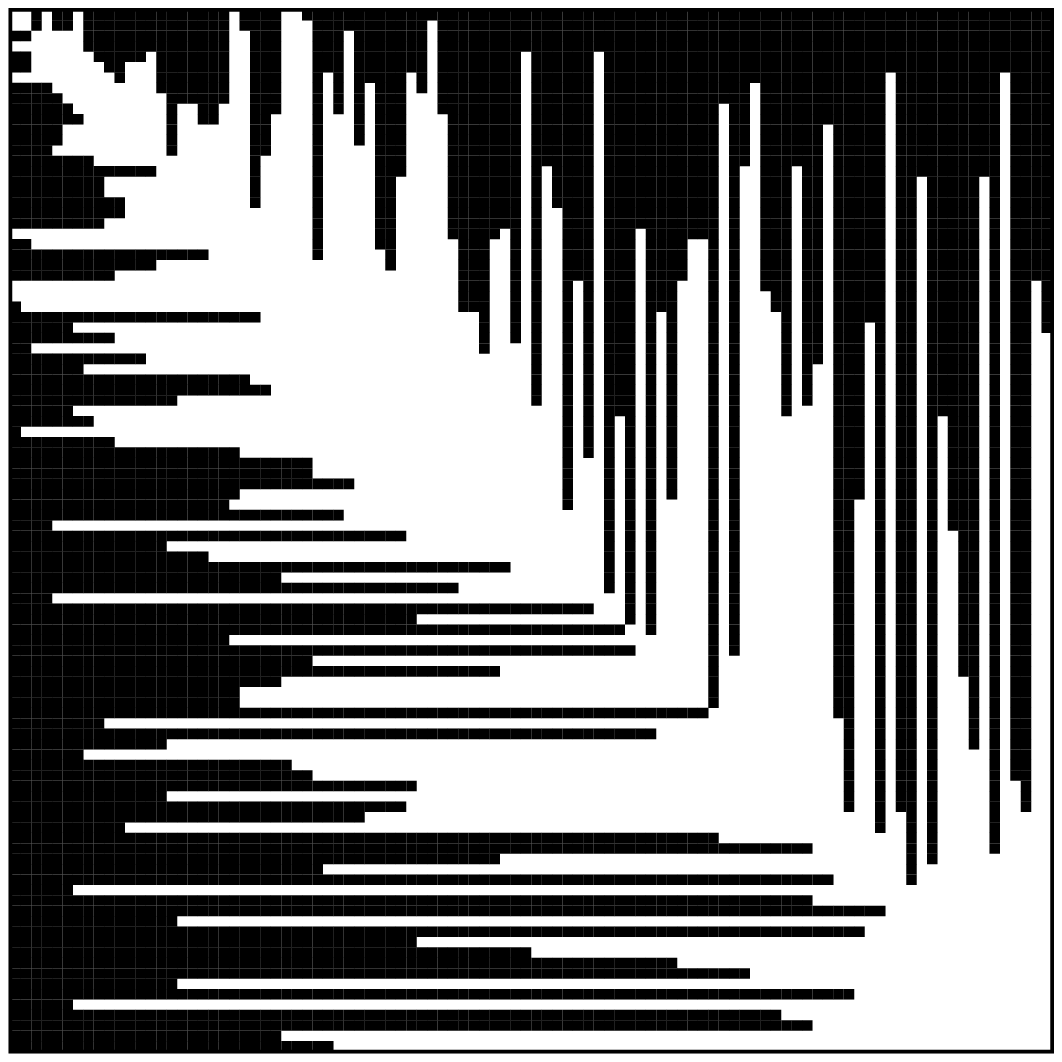}\qquad
  \includegraphics*[height=120pt]{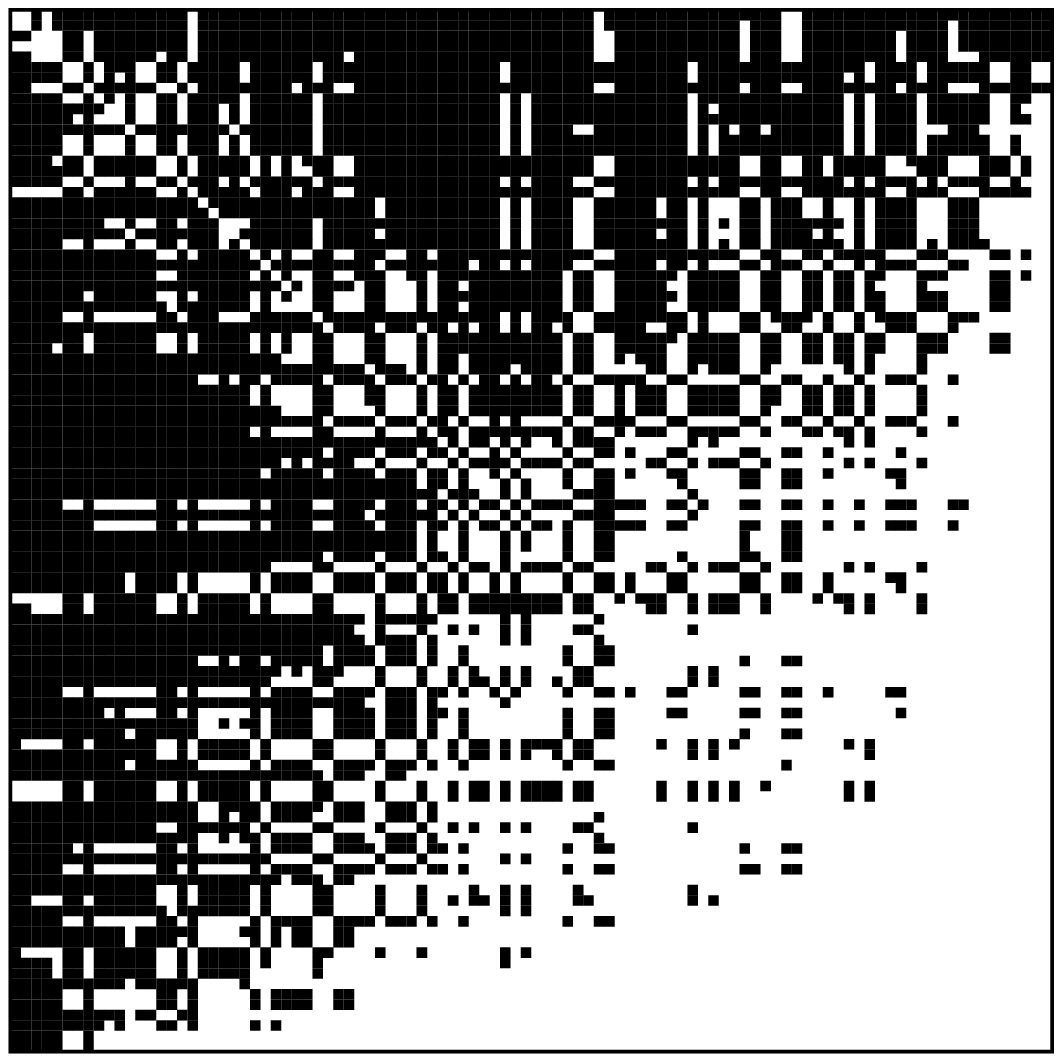}
  \caption{A randomly grown prefix attachment graph with 100 nodes,
  and the same graph with nodes ordered by their degrees.}
  \label{FIG:PREFIX}
\end{figure}

In this construction, it will be more convenient to label the nodes
starting with $1$. At the $n$-th iteration, a new node $n$ is born, a
node $z$ is selected at random, and node $n$ is connected to nodes
$1,\dots,z-1$. We denote the $n$-th graph in the sequence by
$G_n^{\rm pfx}$, and call this graph sequence a {\it randomly grown
prefix attachment graph sequence} (Figure \ref{FIG:PREFIX}).

Again we start with some simple calculations. The probability that
nodes $i<j$ are connected is $\frac{j-i}{j}$ (but these events are
not independent in this case!). The expected degree of $j$ is
therefore
\[
\sum_{i=1}^{j-1} \frac{j-i}{j} + \sum_{i=j+1}^{n} \frac{i-j}{i}
=n-\frac{j}{2}+j\ln\frac{n}{j}+o(n).
\]
The expected number of edges is $n(n-1)/4$.

Looking at the picture, it seems that it tends to some function,
which we can try to figure out similarly as in the case of uniform
attachment graphs. The probability that $i$ and $j$ are connected can
be written in a symmetric form as
\[
\frac{|j-i|}{\max(i,j)}.
\]
If $i=xn$ and $j=yn$, then this is
\[
\frac{|x-y|}{\max(x,y)}.
\]

Does this mean that the function $U(x,y)=|x-y|/\max(x,y)$ is the
limit? Somewhat surprisingly, the answer is negative, which we can
see by computing triangle densities. The probability that nodes
$i<j<k$ form a triangle is
$\bigl(1-\frac{j}{k}\bigr)\bigl(1-\frac{i}{j}\bigr)$ (since if $k$ is
connected to $j$, then it is also connected to $i$). Hence the
expected number of triangles is
\[
\sum_{i<j<k}\left(1-\frac{j}{k}\right)\left(1-\frac{i}{j}\right)
\killtext{&= \sum_{j<k}\left(1-\frac{j}{k}\right)\frac{j-1}{2} =
\sum_{j<k}\left(\frac{j-1}{2}-\frac{j(j-1)}{2k}\right) \\
&=\sum_k\left(\frac{(k-1)(k-2)}{4} - \frac{(k-1)(k-2)}{6}\right)
}
= \frac{1}{6}\binom{n}{3}.
\]
Hence
\[
t(K_3,G_n) = \frac{1}{n^3}\binom{n}{3}\to \frac16.
\]
On the other hand,
\[
t(K_3,U)=\int_{[0,1]^3} \frac{|x-y|}{\max(x,y)}\cdot
\frac{|x-z|}{\max(x,z)} \cdot\frac{|y-z|}{\max(y,z)}\,dx\,dy\,dz.
\]
Since the integrand is independent of the order of the variables, we
can compute this easily:
\[
t(K_3,U)= 6 \int_{0\le x<y<z\le 1}
\left(1-\frac{x}{y}\right)\left(1-\frac{x}{z}\right)
\left(1-\frac{y}{z}\right)\,dx\,dy\,dz=\frac{5}{36}.
\]
So $U$ is not the limit of the sequence $G_n^{\rm pfx}$. On the other
hand, it is not hard to verify that
\begin{equation}\label{WEAKSTAR}
\int_{S\times T} (W_{G_n^{\rm pfx}}-W) \to 0
\end{equation}
for every $S,T\subseteq[0,1]$. Indeed, it is enough to prove this for
sets $S,T$ from a generating set of the $\sigma$-algebra of Borel
sets, e.g. rational intervals. Since there is only a countable number
of these intervals, it suffices to prove that \eqref{WEAKSTAR} holds
with probability $1$ for each fixed $S$ and $T$. It is also easy to
see that it suffices to consider the case $S=T$. For a node $j$ with
$j/n\in S$, let $X_{n,j}$ denote the number of edges $ij$ ($i<j$) in
$G_n^{\rm pfx}$ with $i/n\in S$, and let $X_n=\sum_{j/n\in S}
X_{n,j}$. Then direct computation shows that
\[
\frac1{n^2}\E(X_n)\to \int_{S\times S} U.
\]
Furthermore, the variables $X_{n,j}$ are independent for fixed $n$,
hence the Chernoff--Hoeffding Inequality implies that
$\Pr(|X_n-\E(X_n)|>\eps n^2)$ drops exponentially with $n$. Hence it
follows that $X_n/n^2\to \int_{S\times S} U$ with probability $1$.

So $W_{G_n^{\rm pfx}}\to W$ in the weak-star topology of
$L_\infty[0,1]^2$, but not in our sense. This example also shows that
had we defined convergence of a graph sequence by this convergence in
weak-star topology (after appropriate relabeling), the limit would
not be unique.

Perhaps ordering the nodes by degrees helps? The second pixel picture
in Figure \ref{FIG:PREFIX} suggests that after this reordering, the
functions $W_{G_n^{\rm pfx}}$ converge to some other continuous
function. But again this convergence is only in the weak-star
topology, not in the $\delta_\square$ distance. We'll see that no
continuous function can represent the ``right'' limit: the limit
graphon is $0$-$1$ valued, and it is uniquely determined up to
measure preserving transformations by Theorem \ref{THM:LIMIT}, which
do not change this property.

Is this graph sequence convergent at all? Our computation of the
triangle densities above can be extended to computing the density of
any subgraph, and it follows that the sequence of densities
$t(F,G_n^{\rm pfx})$ is convergent for every $n$. How to figure out
the limit?

Let us label a node born in step $k$, connected to $\{1,\dots,m\}$,
by $(k/n, m/k)\in [0,1]\times[0,1]$. Then we can observe that {\it
nodes with label $(x_1, y_1)$ and $(x_2, y_2)$ are connected if and
only if either $x_1< x_2y_2$ or $x_2< x_1y_1$.}

Consider the function $W:~[0,1]^2\times [0,1]^2\to[0,1]$, given by
\[
W^{\rm pfx}((x_1, y_1),(x_2, y_2))=
  \begin{cases}
    1, & \text{if $x_1< x_2y_2$ or $x_2< x_1y_1$}, \\
    0, & \text{otherwise}.
  \end{cases}
\]

\begin{prop}\label{PROP:PREFIX}
The prefix attachment graphs $G_n^{\rm pfx}$ tend to $W^{\rm pfx}$
with probability $1$.
\end{prop}

\begin{proof}
Let $S_n$ be the (random) set of points in $[0,1]^2$ of the form
$(i/n,z_i/i)$ where $i=1,\dots,n$ and $z_i$ is a uniformly chosen
random integer in $\{1,\dots,i\}$. Then $G_n^{\rm pfx}=\Ge(S_n,W^{\rm
pfx})=\Ha(S_n,W^{\rm pfx})$.

Furthermore, with probability $1$, the sets $S_n$ are well
distributed in $[0,1]^2$. Indeed, for $m\ge 1$, let $J_{m,k}$ denote
the interval $(k/m, (k+1)/m]$, and let $\PP_m$ denote the partition
of $[0,1]^2$ into the sets $J_{m,k}\times J_{m,l}$
($k,l=0,\dots,m-1$). We want to prove that for every fixed $m$ and
$0\le k,l\le m-1$, $|S_n\cap(J_{m,k}\times J_{m,l})|/n\to 1/m^2$ as
$n\to\infty$ with probability $1$. Let
\[
X_i=
   \begin{cases}
    1, & \text{if $(i,z_i)\in J_{m,k}\times J_{m,l}$}, \\
    0, & \text{otherwise},
  \end{cases}
\]
then
\[
|S_n\cap(J_{m,k}\times J_{m,l})| =\sum_{i=1}^n X_i.
\]
We have
\[
\E(X_i)=
  \begin{cases}
     \displaystyle \frac1i\left(\Bigl\lfloor\frac{(l+1)i}{m}\Bigr\rfloor
     -\Bigl\lfloor\frac{li}{m}\Bigr\rfloor\right), &
     \text{if $\displaystyle \frac{k}m\le \frac{i}n\le \frac{k+1}m$},\\
     0, & \text{otherwise},
  \end{cases}
\]
and hence
\begin{align*}
\E|S_n\cap(J_{m,k}\times J_{m,l})| &=\sum_{i\in nJ_{m,k}}
\frac1i\left(\left\lfloor\frac{(l+1)i}{m}\right\rfloor-
\left\lfloor\frac{li}{m}\right\rfloor\right)=\sum_{i\in nJ_{m,k}} \frac1m +O(\log n)\\
&=\frac1m\left(\left\lfloor\frac{(k+1)n}{m}\right\rfloor-
\left\lfloor\frac{kn}{m}\right\rfloor\right) +O(\log
n)=\frac{n}{m^2}+O(\log n).
\end{align*}
Thus
\[
\E\Bigl(\frac1n|S_n\cap(J_{m,k}\times J_{m,l})|\Bigr)
\to\frac1{m^2}\qquad (n\to\infty).
\]
The fact that $|S_n\cap(J_{m,k}\times J_{m,l})|/n \to1/m^2$ with
probability $1$ (not just in expectation) follows by the Law of Large
Numbers, since the $X_i$ are independent.

Thus Lemma \ref{LEM:CONV} applies and proves the Proposition.
\end{proof}

Lemma \ref{LEM:CONV} in fact implies (since $W^{\rm pfx}$ is 0-1
valued) that $W_{G_n^{\rm pfx}}$ tend to $W^{\rm pfx}$ with
probability $1$ in the edit distance, not just in the cut distance.
This means that while the graphs $G_n^{\rm pfx}$ are random, they are
very highly concentrated: two instances of $G_n^{\rm pfx}$ only
differ in $o(n^2)$ edges if overlayed properly (not in the original
ordering of the nodes!). Informally, they have a relatively small
amount of randomness in them, which disappears as $n\to\infty$.
Indeed, $G_n^{\rm pfx}$ is generated using only $O(n\log n)$ bits, as
opposed to, say, $\Ge(n,1/2)$, which is generated using
$\binom{n}{2}$ bits. It would be interesting to explore this
phenomenon.

\begin{figure}[htb]
  \centering
  \includegraphics*[height=120pt,width=120pt]{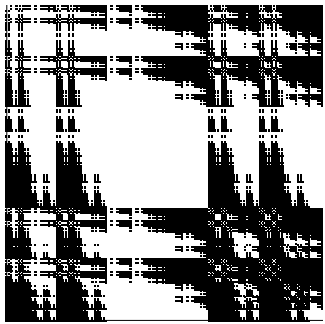}
  \caption{The limit of randomly grown prefix attachment graphs (as a function
on $[0,1]^2$)}\label{FIG:PREFIX-LIM}
\end{figure}

Proposition \ref{PROP:PREFIX} gives a nice and simple representation
of the limit object with the underlying probability space $[0,1]^2$
(with the uniform measure). If we want a representation on $[0,1]$,
we can map $[0,1]$ into $[0,1]^2$ by a measure preserving map $\phi$;
then $W_{\rm pfx}^\phi(x,y)=W^{\rm pfx}(\phi(x),\phi(y))$ gives a
representation of the same graphon as a 2-variable function. For
example, using the map $\phi$ that separates even and odd bits of
$x$, we get the fractal-like picture in Figure \ref{FIG:PREFIX-LIM}.

It is interesting to note that the graphs $\Ge(n,W)$ form another
(different) sequence of random graphs tending to the same limit $W$
with probability 1.

\subsection{Preferential attachment graph on n fixed nodes}\label{SEC:PFR-ATT}

A {\it preferential attachment graph with $n$ fixed nodes and $m$
edges} $\PAG(n,m)$ is the random graph obtained by the following
procedure. Let $v_1\dots v_n$ be a set of nodes. We extend this
sequence one by one by picking an element of the current sequence
randomly and uniformly, and append a copy of it at the end. We repeat
this until $2m$ further elements have been added. So we get a
sequence $v_1\dots v_nv_{n+1}\dots v_{n+2m}$.

Now we connect nodes $v_{n+2k-1}$ and $v_{n+2k}$ for $k=1,2,\dots,m$,
to get $\PAG(n,m)$. (Note that $\PAG(n,m)$ may have multiple edges
and loops, which we have to live with for the time being).

Another way of describing this construction is to view it as adding
edges one by one, where the probability of adding an edge connecting
$u$ and $v$ is proportional to the product of their degrees. To be
more precise, the probability that the $(k+1)$-st edge connects $u$
and $v$ is
\[
  \begin{cases}
    \displaystyle\frac{2(d_k(u)+1)(d_k(v)+1)}{(n+2k)(n+2k+1)}
     & \text{if $u\not= v$}, \\[12pt]
    \displaystyle\frac{(d_k(u)+1)(d_k(u)+2)}{(n+2k)(n+2k+1)}
     & \text{if $u=v$},
  \end{cases}
\]
where $d_k(u)$ is the current degree of the node (adding 1 to the
degree is needed to start the procedure at all; adding 2 to the
second factor in the case when $u=v$ is a minor trick that makes
everything come out nicer).

\begin{figure}[htb]
  \centering
  \includegraphics*[height=120pt,width=120pt]{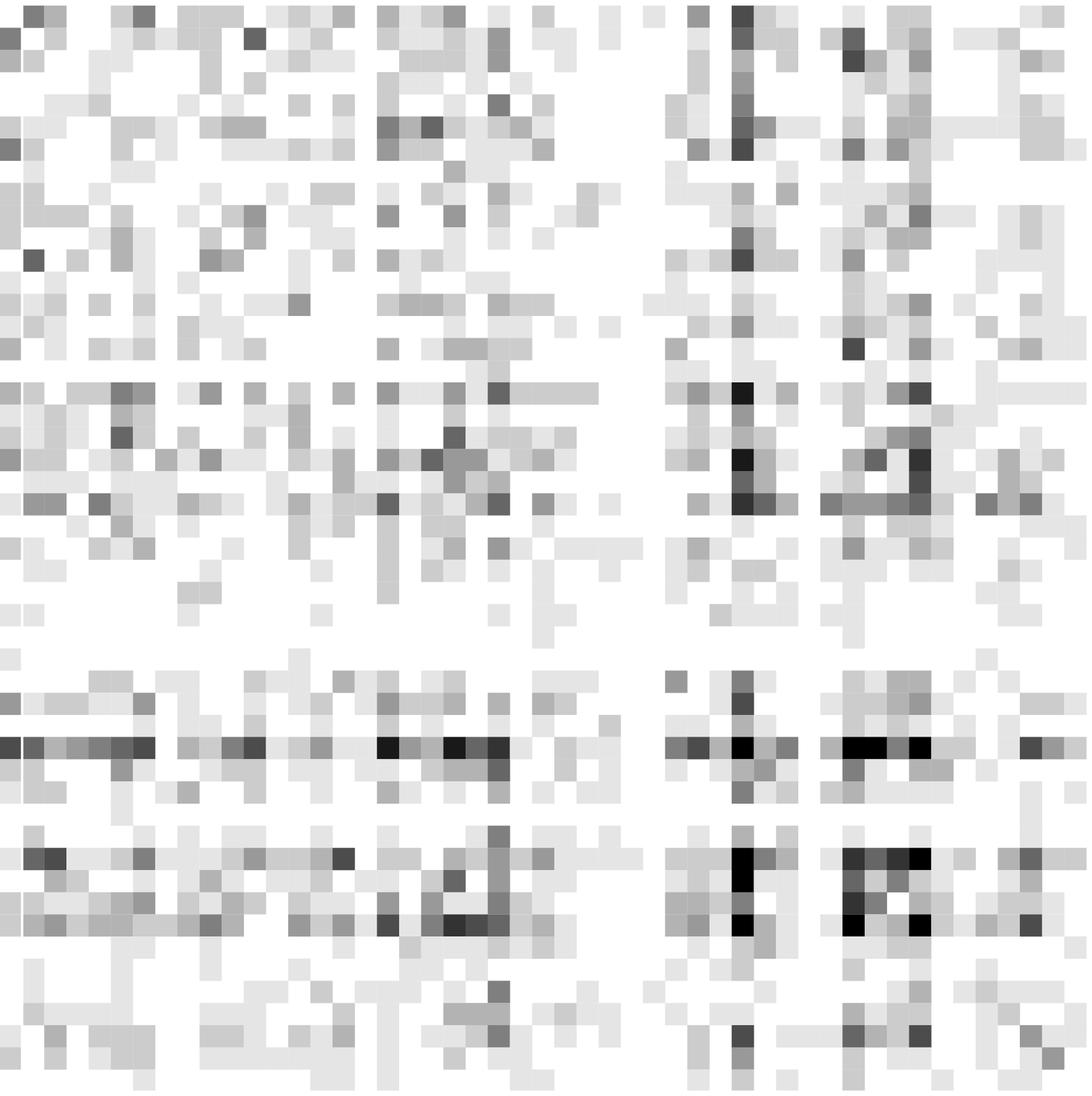}\qquad
  \includegraphics*[height=120pt,width=120pt]{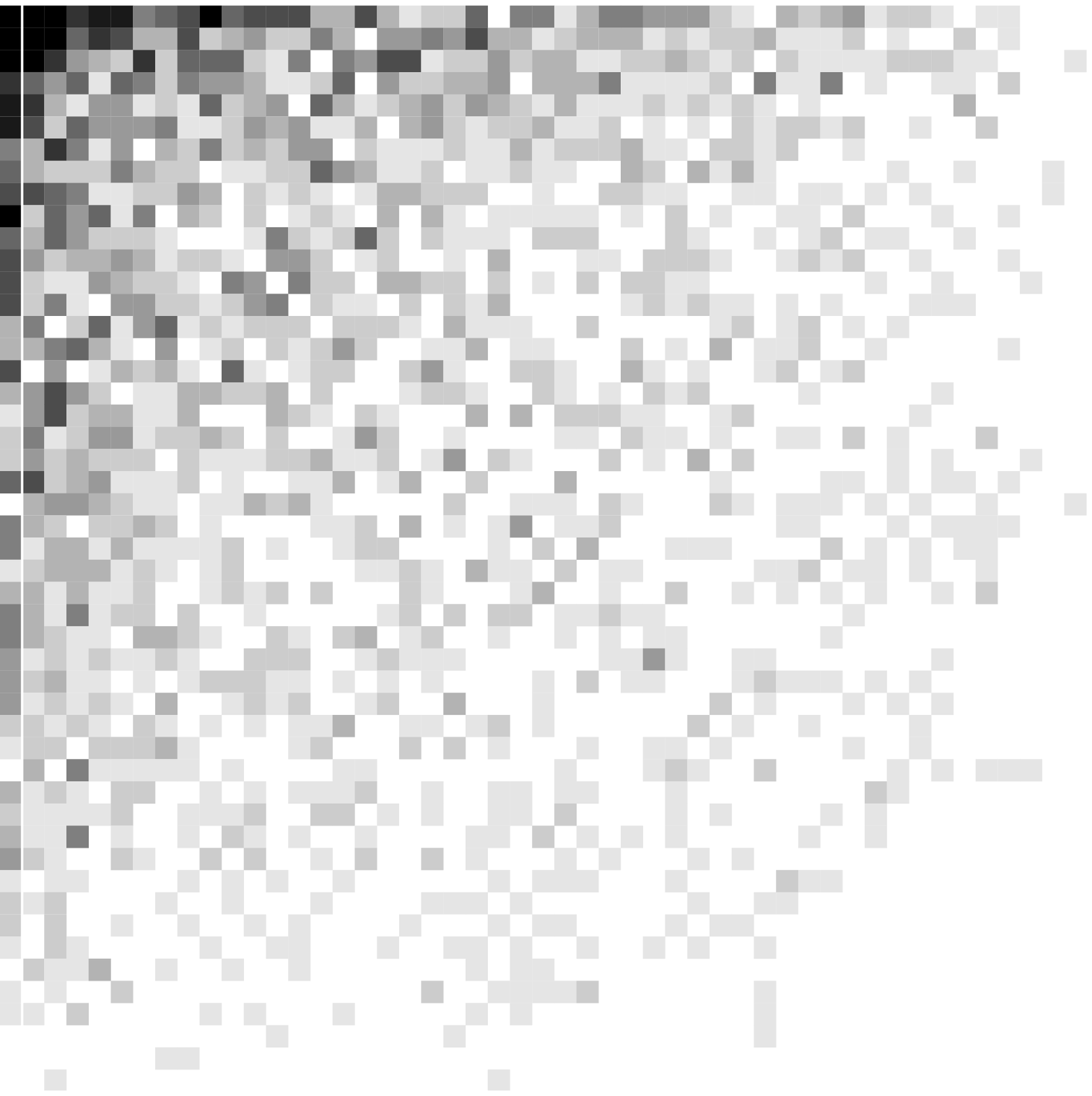}
  \caption{(a) A preferential attachment graph $\PAG(50,1000)$.
  Darkness of a pixel indicates multiplicity of
  the edge. (b) The same graph with the nodes ordered by decreasing
  degrees.}
  \label{FIG:PREF}
\end{figure}

Preferential attachment graphs are motivated by the (sparse)
Albert--Barab\'asi graphs \cite{AlBa}, and they have been studied in
detail by Pittel \cite{Pit}.

The somewhat awkward definition of preferential attachment graphs is
justified by the following nice properties. First, let us compute the
probability that this process yields a multigraph $G$ on $V(G)=[n]$,
with degrees $d_1,\dots,d_n$, with $m$ edges and $m'$ non-loop edges.
Fix any order of the edges, and for the non-loop edges fix an order
in which their endpoints were inserted (i.e., an orientation of $G$).
Then the probability that $G$ arises this way is
\begin{equation}\label{PROB-G-ORD}
\frac{d_1!\dots d_n!}{n(n+1)\dots (n+2m-1)}.
\end{equation}
Summing over all orientations and orderings of the edges, we get that
the probability that $\PAG(n,m)=G$ is
\begin{equation}\label{PROB-G}
m! 2^{m'} \frac{d_1!\dots d_n!}{n(n+1)\dots (n+2m-1)}.
\end{equation}

An important observation we can make from this computation is the
following:

\begin{lemma}\label{INVARIANT}
Conditioning on the graph $G(n,m)$, all the $2^{m'}m!$ possibilities
in which the edges could have been inserted have the same
probability.
\end{lemma}

We can use this lemma to determine the expected subgraph densities in
$\PAG(n,m)$. For two multigraphs $F$ and $G$, let $inj(F,G)$ denote
the number of embeddings of $f$ into $G$, i.e., the number of pairs
$(\phi,\psi)$ of injective maps $\phi:~V(F)\to V(G)$ and
$\psi:~E(F)\to E(G)$ that preserve incidence. Let
\[
t_\inj(F,G)=\frac{\inj(F,G)}{(n)_k},
\]
where $k=|V(F)|$ and $n=|V(G)|$.

Let $F$ be a multigraph on $V(G)=[k]$, with degrees $r_1,\dots,r_k$,
with $l$ edges and $l'$ non-loop edges. Fix an order of the edges of
$F$ and also an orientation $\sigma$ of the non-loop edges of $F$ as
above. Let $\overrightarrow{e}_1, \dots,\overrightarrow{e}_m$ be the
order and orientation in which $\PAG(n,m)$ arises. Let
$p(\sigma,v_1,\dots,v_k, j_1,\dots,j_l)$ denote the probability that
edges $\overrightarrow{e}_{j_1}, \dots,\overrightarrow{e}_{j_l}$ form
a copy of $F$ on nodes $v_1,\dots,v_k$ (with the given labeling of
the nodes, the given order of the edges, and the given orientation).
By Lemma \ref{INVARIANT}, this number is the same for any $l$-tuple
$(j_1,\dots,j_l)$, and trivially, it is the same for every $k$-tuple
$(v_1,\dots,v_k)$. Hence
\begin{align*}
\E\bigl(\inj(F,\PAG(n,m))\bigr)&= \sum_{v_1,\dots,v_k}
\sum_{j_1,\dots,j_l}\sum_{\sigma} p(\sigma,v_1,\dots,v_k,
j_1,\dots,j_l)\\
&= (n)_k(m)_l 2^{l'} p(\sigma_0,1,\dots,k,1,\dots,l),
\end{align*}
where $\sigma_0$ is any fixed orientation of $F$. By
(\ref{PROB-G-ORD}), we have
\[
p(\sigma_0,1,\dots,k,1,\dots,l)=\frac{r_1!\dots r_k!}{(n+2l-1)}_{2l},
\]
and so
\begin{equation}\label{PAG-T}
\E\bigl(t_\inj(F,\PAG(n,m))\bigr)=\frac{1}{(n)_k}(n)_k(m)_l
2^{l'}\frac{r_1!\dots r_k!}{(n+2l-1)}_{2l}=2^{l'}r_1!\dots r_k!
\frac{(m)_l}{(n)_{2l}}.
\end{equation}

Suppose that $n,m\to\infty$ so that $m\sim cn^2/2$. Then
\[
\E\bigl(t_\inj(F,\PAG(n,m))\bigr)\sim 2^{l'}r_1!\dots r_k!
\frac{m^l}{n^{2l}}\longrightarrow 2^{l'-l}c^lr_1!\dots r_k!\,.
\]
If we assume that $F$ has no loops, then
\[
\E\bigl(t_\inj(F,\PAG(n,m))\bigr)\longrightarrow c^lr_1!\dots r_k!\,.
\]
Using high concentration results, one can show that this convergence
holds not only in expectation, but with probability 1,
\[
t_\inj(F,\PAG(n,m))\longrightarrow c^lr_1!\dots r_k!\,.
\]

Note that the relation $t_\inj(F,\PAG(n,m))\sim t(F,\PAG(n,m))$ does
not hold in general if $F$ has multiple edges. In fact, it is easy to
see that
\[
t_\inj(F,\PAG(n,m))\sim \sum_{F'} t(F,\PAG(n,m))\prod_{i,j\in V(F)}
m'_ij! \Bigl\{{m_{ij}\atop m'_{ij}}\Bigr\},
\]
where $F'$ ranges through all multigraphs obtained from $F$ by
reducing the edge multiplicities (not strictly, but keeping at least
one copy of each edge), $m_{ij}$ and $m_{ij}'$ denote the
multiplicities of the edge $ij$ in $F$ and $F'$, respectively, and
$\{ {a\atop b}\}$ denotes the Stirling number of the second kind. For
example, if $K_2^{(2)}$ denotes the graph on two nodes having two
parallel edges, then
\[
t(K_2^{(2)},\PAG(n,m))\sim t_\inj(K_2^{(2)},\PAG(n,m))
+t_\inj(K_2,\PAG(n,m)).
\]

Let $L_c(x,y)=c(\ln x)(\ln y)$. Then for a multigraph $F$ without
loops we have
\[
t(F,L_c)=\int_{[0,1]^k} \prod_{ij\in E(F)} W(x_i,x_j)\,dx =
\int_{[0,1]^k} c^l\prod_{i=1}^k (\ln x_i)^{r_i}\,dx =c^lr_1!\dots
r_k!\,.
\]

This implies that the limit of preferential attachment graphs
$\PAG(n, cn^2)$, with probability 1, is described by the function
$L_c$. To be precise, the graphs $\PAG(n, cn^2)$ have multiple edges,
and hence the theory of convergent graph sequences developed in
\cite{BCLSV1,BCLSV2} does not apply, but the computations above show
that the convergence does occur in at least one possible sense.

\begin{prop}\label{PROP:PREFER}
If $m(n)=(c+o(1))n^2$, then with probability $1$,
$t_\inj(F,\PAG(n,m))\to t(F,L_c)$ for every multigraph $F$ without
loops.
\end{prop}

Let $\SPAG(n,cn^2)$ denote the {\it simplified preferential
attachment graph} obtained from $\PAG(n,cn^2)$ by deleting loops and
keeping only one copy of parallel edges. L.~Szak\'acs \cite{Szak}
proved that this sequence of graphs is convergent with probability
$1$, and its limit is the function $1-\exp(-c\ln x\ln y)$.

\section{Convergence to a prescribed function}

Lemma \ref{LEM:WRAND} gives a way to construct a randomly growing
graph sequence converging to a given function $W$. Let
$s_1,s_2,\dots\in\Omega$ be independent random samples from $\pi$,
and let $S_n=\{s_1,\dots,s_n\}$. We can construct $\Ge(S_n,W)$ by
taking $G(S_{n-1},W)$, adding $s_n$ as a new node, and connecting
$s_n$ to $s_i$ with probability $W(s_n,s_i)$. Then
$\Ge(S_1,W),\Ge(S_2,W),\dots$ is a randomly growing sequence of
graphs, and by Lemma \ref{LEM:WRAND}, we have $\Ge(S_n,W)\to W$ with
probability $1$.

\killtext{ ??? Another way of describing how $G_n=\Ge(S_n,W)$ is
obtained from the graph $G_{n-1}=\Ge(S_{n-1},W)$ is the following.
For $J\subseteq V(G_{n-1})$, let $\bar{G}_{n-1,J}$ be obtained from
$\bar{G}_{n-1}$ by adding a new node and connecting it to the points
in $J$ by red edges, and to the other nodes of $V(G_{n-1})$ by blue
nodes. Then
\[
p(J)=\frac{t(\bar{G}_{n-1,J},W)}{t(\bar{G}_{n-1},W)}
\]
defines a probability distribution on $2^{V(G_{n-1})}$, and $G_n$ can
be obtained from $G_{n-1}$ by adding a new node and connecting it to
a set of old nodes $J$ with probability}

However, one can have several objections to this method: First, the
new edges are not added independently of each other. Second, even if
$\Omega=[0,1]$, and the function $W$ is, say, continuous and
monotone, the order in which the nodes of $\Ge(S_n,W)$ are born is
random, and has nothing to do with the order of the points
$s_i\in[0,1]$ representing them. In other words, to get a labeling
for which $W_{\Ge(S_n,W)}\to W$ in the $\|.\|_\square$ norm, we have
to reorder the nodes.

It may be interesting to consider rules for generating randomly
growing graph sequences $(G_n)$ with a prescribed limit function $W$
for which these objections cannot be raised. Given a function
$W\in\WW_0$, monotone decreasing in each variable, construct a
randomly growing simple graph sequence $(G_1,G_2,\dots)$ as follows.
$G_1$ is a single node labeled $1$. For $n>1$, define
\[
p_{n,j}=W(\tfrac{j}{n},1),\qquad p_{n,ij}=
\frac{W(\frac{i}{n},\frac{j}{n})-W(\frac{i}{n-1},
\frac{j}{n-1})}{1-W(\frac{i}{n-1},\frac{j}{n-1})}.
\]
To get $G_n$ from $G_{n-1}$, we add a new node $n$, connect it to
each node $j<n$ with probability $p_{n,j}$, and connect any two
nonadjacent nodes $i,j<n$ with probability $p_{n,ij}$. All these
decisions are independent. The monotonicity of $W$ implies that $0\le
p_{n,ij}\le 1$ is a legal probability.

\begin{prop}\label{PROP:}
The sequence of graphs $G_n$ constructed above has the property that
$W_{G_n}\to W$ in the $\|.\|_\square$ norm.
\end{prop}

\begin{proof}
The probability that nodes $i<j$ are not connected in $G_n$ is
\begin{align*}
(1-p_{j,i})(1-p_{j+1,ij})\cdots(1-p_{n,ij})
&=\bigl(1-W(\tfrac{i}{j},1)\bigr)
\frac{1-W(\frac{i}{j+1},\frac{j}{j+1})}{1-W(\frac{i}{j},\frac{j}{j})}
\cdots
\frac{1-W(\frac{i}{n},\frac{j}{n})}{1-W(\frac{i}{n-1},\frac{j}{n-1})}\\
&=1-W(\tfrac{i}{n},\tfrac{j}{n}),
\end{align*}
and hence the probability that they are adjacent is
$W(\tfrac{i}{n},\tfrac{j}{n})$. Thus $G_n$ is the graph $\Ge(S_n,W)$,
where $S_n=\{\frac1n,\tfrac2n,\dots\tfrac{n-1}n\}$. It is trivial
that this sequence of sets is well distributed in $[0,1]$, and since
$W$ is almost everywhere continuous, it follows by Lemma
\ref{LEM:CONV} that $G_n\to W$ with probability $1$.
\end{proof}

The convergent sequences discussed in Sections \ref{SEC:UNIF-ATT} and
\ref{SEC:R-ATT} are special cases of this construction. A more
general nice case is when $W=1-U$, where $U$ is homogeneous of some
degree: $U(\lambda x, \lambda y)= \lambda^c U(x,y)$ with some $c\ge
0$. When a new node $n$ is born we connect it to node $i<n$ with
probability $W(\frac in,1)$, and then at each further step, we
connect any two nonadjacent nodes with probability
$1-\bigl(\frac{n-1}{n}\bigr)^c$.

\end{document}